\documentclass[12pt]{article}
\usepackage{graphicx}
\usepackage{amsmath,amsthm,amssymb,enumerate}%, esint}
\usepackage{euscript,mathrsfs}
\usepackage{dsfont}
\usepackage[left=2cm,right=2cm,top=3.5cm,bottom=3.5cm]{geometry}
\usepackage{color}
\catcode`\@=11 \@addtoreset{equation}{section}

\catcode`\@=12

\allowdisplaybreaks

\newtheorem{Theorem}{Theorem}[section]
\newtheorem{Proposition}[Theorem]{Proposition}
\newtheorem{Lemma}[Theorem]{Lemma}
\newtheorem{Corollary}[Theorem]{Corollary}

\theoremstyle{definition}
\newtheorem{Definition}[Theorem]{Definition}

\newtheorem{Remark}[Theorem]{Remark}

\newcommand{\bTheorem}[1]{
\begin{Theorem} \label{T#1} }
\newcommand{\eT}{\end{Theorem}}

\newcommand{\bProposition}[1]{
\begin{Proposition} \label{P#1}}
\newcommand{\eP}{\end{Proposition}}

\newcommand{\bLemma}[1]{
\begin{Lemma} \label{L#1} }
\newcommand{\eL}{\end{Lemma}}

\newcommand{\bCorollary}[1]{
\begin{Corollary} \label{C#1} }
\newcommand{\eC}{\end{Corollary}}

\newcommand{\bRemark}[1]{
\begin{Remark} \label{R#1} }
\newcommand{\eR}{\end{Remark}}

\newcommand{\bDefinition}[1]{
\begin{Definition} \label{D#1} }
\newcommand{\eD}{\end{Definition}}

\newcommand{\Nu}{\mathcal{V}_{t,x}}

\newcommand{\bfphi}{\boldsymbol{\varphi}}

\newcommand{\bFormula}[1]{
\begin{equation} \label{#1}}
\newcommand{\eF}{\end{equation}}

\newcommand{\Ov}[1]{\overline{#1}}

\newcommand{\vr}{\varrho}

\newcommand{\tvr}{\tilde \vr}
\newcommand{\tvu}{{\tilde \vu}}
\newcommand{\tvt}{\tilde \vt}

\newcommand{\vt}{\vartheta}
\newcommand{\vu}{\vc{u}}
\newcommand{\vm}{\vc{m}}

\newcommand{\vM}{\vc{M}}

\newcommand{\vn}{\vc{n}}

\newcommand{\vc}[1]{{\bf #1}}

\newcommand{\Div}{{\rm div}_x}
\newcommand{\Grad}{\nabla_x}

\newcommand{\dx}{\,{\rm d} {x}}

\newcommand{\dt}{\,{\rm d} t }

\newcommand{\vU}{\vc{U}}

\newcommand{\intO}[1]{\int_{\Omega} #1 \ \dx}

\newcommand{\ep}{\varepsilon}

\definecolor{Cgrey}{rgb}{0.85,0.85,0.85}
\definecolor{Cblue}{rgb}{0.50,0.85,0.85}
\definecolor{Cred}{rgb}{1,0,0}
\definecolor{fancy}{rgb}{0.10,0.85,0.10}

\newcommand\Cbox[2]{%
    \newbox\contentbox%
    \newbox\bkgdbox%
    \setbox\contentbox\hbox to \hsize{%
        \vtop{
            \kern\columnsep
            \hbox to \hsize{%
                \kern\columnsep%
                \advance\hsize by -2\columnsep%
                \setlength{\textwidth}{\hsize}%
                \vbox{
                    \parskip=\baselineskip
                    \parindent=0bp
                    #2
                }%
                \kern\columnsep%
            }%
            \kern\columnsep%
        }%
    }%
    \setbox\bkgdbox\vbox{
        \color{#1}
        \hrule width  \wd\contentbox %
               height \ht\contentbox %
               depth  \dp\contentbox
        \color{black}
    }%
    \wd\bkgdbox=0bp%
    \vbox{\hbox to \hsize{\box\bkgdbox\box\contentbox}}%
    \vskip\baselineskip%
}

\date{}

\begin{document}

%%%%%%%%%%%%%%%%%%%%%%%%%%%%%%%%

\title{On a singular limit for stratified compressible fluids}

\author{Gabriele Bruell $\phantom{}^*$ \and Eduard Feireisl
\thanks{The research leading to these results has received funding from the
European Research Council under the European Union's Seventh
Framework Programme (FP7/2007-2013)/ ERC Grant Agreement
320078. The Institute of Mathematics of the Academy of Sciences of
the Czech Republic is supported by RVO:67985840.}
}

\date{\today}

\maketitle

\bigskip

\centerline{Department of Mathematical Sciences, Norwegian University of Science and Technology}

\centerline{Alfred Getz vei 1, 7491 Trondheim, Norway}

\bigskip

\centerline{Institute of Mathematics of the Academy of Sciences of the Czech Republic}

\centerline{\v Zitn\' a 25, CZ-115 67 Praha 1, Czech Republic}

\bigskip

\begin{abstract}

We consider a singular limit problem for the complete compressible Euler system in the low Mach and strong stratification regime. We identify the limit problem - the anelastic Euler system - in the case of well prepared initial data. The result holds in the large class of the dissipative measure--valued solutions of the primitive system.
Applications are discussed to the driven shallow water equations.

\end{abstract}


%\tableofcontents

\allowdisplaybreaks

\section{Introduction}
\label{i}

The following system of equations arises in a number of real world applications, in particular in certain
astrophysical and meteorological models (see e.g. the survey by
Klein \cite{Klein}):
\begin{equation} \label{i1}
\partial_t \vr + \Div (\vr \vu) = 0,
\end{equation}
\begin{equation} \label{i2}
\partial_t (\vr \vu) + \Div (\vr \vu \otimes \vu) + \frac{1}{\ep^2} \Grad p(\vr, \vt) = \frac{1}{\ep^2} \vr \Grad F,
\end{equation}
\begin{equation} \label{i3}
\partial_t \left( \frac{1}{2} \vr |\vu|^2 + \frac{1}{\ep^2} \vr e(\vr, \vt) \right) + \Div \left[ \left( \frac{1}{2} \vr |\vu|^2 + \frac{1}{\ep^2} \vr e(\vr, \vt) \right) \vu \right] +
\frac{1}{\ep^2} \Div (p(\vr, \vt) \vu ) = \frac{1}{\ep^2} \vr \Grad F \cdot \vu.
\end{equation}
The equations (\ref{i1}), (\ref{i2}), and (\ref{i3}) represent a mathematical formulation of the conservation of mass, momentum, and energy, respectively, of a compressible inviscid fluid driven by a potential force $\Grad F$. Here, the state of the fluid at a time $t$ and a spatial position $x$ is given by its mass density $\vr = \vr(t,x)$, the macroscopic
velocity $\vu = \vu(t,x)$ and the (absolute) temperature $\vt = \vt(t,x)$. The pressure $p = p(\vr, \vt)$ and the internal energy density $e = e(\vr, \vt)$ are given explicitly through an equation of state.

To close the system we specify the physical domain - an infinite slab, ``periodic'' in the horizontal variable:
\[
\Omega = \left( [-1,1]|_{\{ -1,1 \}} \right)^2 \times [0,1],
\]
supplemented with the impermeability condition
\begin{equation} \label{i4}
\vu \cdot \vc{n}|_{\partial \Omega} = 0, \ \mbox{meaning},\ u^3(t,x_1,x_2,0) = u^3(t, x_1, x_2, 1) = 0.
\end{equation}
Later, to extend the range of possible applications of our result, we consider a slightly more general setting, $\Omega \subset R^N$, $N=2,3$,
 - a bounded regular domain - supplemented with the impermeability
condition
\begin{equation} \label{i4a}
\vu \cdot \vc{n}|_{\partial \Omega}= 0.
\end{equation}

Problem (\ref{i2}--\ref{i4}) contains a small positive parameter $\ep > 0$. Our aim is to identify the limit problem for $\ep \to 0$. Rather surprisingly, the limit problem is not unique and depends on the choice of the initial data. To see this, let us examine the associated static system
\begin{equation} \label{i5}
\Grad p(\vr, \vt) = \vr \Grad F.
\end{equation}
To simplify presentation, we suppose that $p$ satisfies the standard Boyle--Mariotte law,
\[
p(\vr, \vt) = \vr \vt,
\]
and that
\[
F(x) = - x_3.
\]
Accordingly, the pressure $p$ in (\ref{i5}) depends only on the vertical variable $x_3$ and problem (\ref{i5}) reduces to
\begin{equation} \label{i6}
\frac{1}{\vr} \partial_{x_3} (\vr \vt) = - 1.
\end{equation}

\subsection{Isothermal limit}

Suppose, that $\vt = \Ov{\vt} > 0$ - a positive constant. Then the stationary problem (\ref{i6}) can be explicitly solved,
\[
\vr = \tvr (x_3),\ \ \tvr(x_3) = c_M \exp \left( - \frac{x_3}{\Ov{\vt}} \right), \ c_M > 0,
\]
where the value of the constant $c_M$ is uniquely determined by prescribing the total mass $M = \intO{ \tvr }$.
The limit system for the isothermal case has been identified in \cite{FeKlKrMa}. It turns out that the limit velocity field
$\vc{U}$ has only two components,
\[
\vc{U}(x_1,x_2,x_3) = \left[ U^1 (x_1,x_2,x_3), U^2(x_1,x_2,x_3),0 \right] \equiv \left[\vc{U}_h (x_1,x_2,x_3) ,0\right],
\]
where, for any fixed $x_3$, the field $\vc{U}_h(\cdot, x_3)$ satisfies the incompressible Euler system
\begin{equation} \label{i7}
{\rm div}_h \vc{U}_h = 0,
\end{equation}
\begin{equation} \label{i8}
\partial_t \vc{U}_h + \vc{U}_h \cdot \nabla_h \vc{U}_h + \nabla_h \Pi = 0,
\end{equation}
in $(0,T) \times \left([-1,1]|_{\{ -1, 1 \}} \right)^2$.
Here the subscript $h$ used with a symbol of differential operator indicates that the latter applies only in the horizontal variable $x_h= [ x_1,x_2]$.

Note that (\ref{i7})--(\ref{i8}) is the 2D incompressible Euler system, parameterized by the vertical variable $x_3$, that admits smooth global in time solutions for any sufficiently regular initial data. The convergence holds provided the initial data are \emph{well prepared}, meaning sufficiently close to the equilibrium state and the
2D-initial velocity:
\[
\vr(0, \cdot) \approx \tvr,\ \vt(0, \cdot) = \Ov{\vt},\ \vu(0,\cdot) \approx [\vc{U}_{0,h}, 0],\ {\rm div}_h \vc{U}_{0,h} = 0.
\]

A heuristic argument why the singular limit exhibits only horizontal motion is based on the entropy equation associated to (\ref{i1}--\ref{i3}). Introducing the caloric EOS,
\[
e(\vr, \vt) = c_v \vt,\ c_v > 0,
\]
with the associated entropy
\[
s(\vr, \vt) = \log \left( \frac{\vt^{c_v}}{\vr} \right),
\]
we deduce from (\ref{i1}--\ref{i3}) the entropy balance equation
\begin{equation} \label{i9}
\partial_t (\vr s(\vr, \vt)) + \Div (\vr s(\vr, \vt) \vu ) = 0
\end{equation}
as long as all quantities in question are smooth enough. As $\vr \approx \tvr$, $\vt \approx \Ov{\vt}$ in the asymptotic limit, the equation of continuity (\ref{i1}) reduces to
\begin{equation} \label{i10}
\Div (\tvr \vc{U} ) = 0,
\end{equation}
while (\ref{i9}) gives rise
\begin{equation} \label{i11}
\Div (\tvr s(\tvr, \Ov{\vt} ) \vc{U} ) = 0.
\end{equation}
It is easy to check that (\ref{i10}), (\ref{i11}) are compatible only if $U^3 \equiv 0$.

\subsection{Isentropic limit}

In the present paper we focus on the isentropic limit, studied in the context of the Navier--Stokes fluid in \cite{FeKlNoZa}. We suppose that the stationary state is isentropic,
specifically, $\vr = \tvr(x_3)$, $\vt = \tvt(x_3)$ such that
\begin{equation} \label{i12}
\frac{ \tvt^{c_v} }{\tvr} = a > 0 \ \mbox{- a positive constant},\quad a^\frac{1}{c_v} \partial_{x_3} \tvr^\gamma = - \tvr ,\quad \gamma = 1 + \frac{1}{c_v}.
\end{equation}
Clearly, (even though not relevant for the subsequent analysis) from \eqref{i12} we can compute the solution to the static problem explicitly,
\begin{equation}
\tvr(x_3) = \left( c_M - \frac{\gamma - 1}{\gamma a^{1/c_v}} x_3 \right)^{\frac{1}{\gamma - 1}},
\end{equation}
where $c_M>0$ is a constant determined uniquely by the total mass.

Our goal will be to show that the limit velocity field $\vc{U}$ is now described by the \emph{anelastic Euler system}:
\begin{equation} \label{i13}
\Div (\tvr \vc{U} ) = 0,
\end{equation}
\begin{equation} \label{i14}
\partial_t \vc{U} + \vc{U} \cdot \Grad \vc{U} + \Grad \Pi = 0
\end{equation}
in $(0,T) \times \Omega$
as long as  the initial data are sufficiently close to the static solution
\[
\vr(0, \cdot) \approx \tvr,\ \vt(0, \cdot) \approx \tvt,\ \mbox{and}\ \vc{u}(0, \cdot) \approx \vc{U}_0,\ \ \Div( \tvr \vc{U}_0 ) = 0.
\]
Note that the system (\ref{i13}), (\ref{i14}) admits smooth solutions for smooth initial data defined on a maximal life span $[0,T_{\rm max})$. Moreover, $T_{\rm max} = \infty$ if $N=2$, see Oliver \cite{Oli}.

Finally, we point out that the \emph{target systems} (\ref{i7})--(\ref{i8}) and (\ref{i13})--(\ref{i14}) describing the asymptotic fluid velocity are apparently different, while the \emph{primitive
system} (\ref{i1}--\ref{i3}) is the same. Thus the convergence result is quite sensitive with respect to the choice of the initial data.

Our goal in this paper is to justify the asymptotic limit in the isentropic case. To this end, we introduce the dissipative measure--valued solutions to the primitive system
(\ref{i1}--\ref{i3}) and review their basic properties in Section \ref{p}. The main result is stated in Section \ref{M}. The proof of the main result, based on the application of the relative energy inequality in the context of measure--valued solutions, is given in Section \ref{PM}. Applications to the shallow water equations and related problems are discussed in Section \ref{A}.

\section{Preliminaries, solutions of the primitive system}
\label{p}

As the primitive system (\ref{i1}--\ref{i3}) is non--linear hyperbolic, the existence of global in time \emph{regular} solutions is in general precluded by formation of singularities as shock waves
in a finite time lap. The existence of weak solutions, on the other hand, is not known for general initial data. In addition, the problem is not well posed for bounded initial data in $L^\infty$, even in the class of entropy admissible solutions, see \cite{FeKlKrMa}.
Similarly to \cite{FeKlKrMa},
we examine the problem in the class of \emph{dissipative measure--valued (DMV) solutions} introduced recently in \cite{BreFei17}. The advantage of such an approach to a singular limit problem is obvious:

\begin{itemize}
\item
DMV solutions to the primitive system exist globally in time.
The result is therefore restricted only by the life-span of the target system that may be finite.

\item The class of solutions of the primite system is very general (large), therefore the result concerning the asymptotic limit is unconditional and in a way the best possible.

\end{itemize}

\subsection{Dissipative measure--valued solutions to the primitive system}

We recall the definition of (DMV) solutions to the Euler system (\ref{i1}--\ref{i3}) introduced in \cite{BreFei17}. A suitable phase space is spanned by the values of
the density $\vr$, the momentum $\vc{m} = \vr \vu$, and the internal energy, or, equivalently, the pressure $p = \vr \vt$. Accordingly, we introduce
\[
\mathcal{F} = \left\{ [\vr, \vc{m}, p] \ \Big| \ \vr \in [0, \infty),\ \vc{m} \in R^N,\ p \in [0, \infty) \right\},\ N=2,3.
\]

A dissipative measure--valued solution of the problem (\ref{i1}--\ref{i4}) is a parameterized family of probability measures
\[
\left\{ \mathcal{V}_{t,x} \right\}_{t \in [0,T]; x \in \Omega},\ \Nu : L^\infty_{\rm weak-(*)} ((0,T) \times \Omega; \mathfrak{P} (\mathcal{F}))
\]
satisfying:
\begin{itemize}
\item
\begin{equation} \label{p1}
\left[ \intO{ \left< \Nu ; \vr \right> \varphi} \right]_{t = 0}^{t = \tau} = \int_0^\tau \intO{ \Big[ \left< \Nu; \vr \right>
\partial_t \varphi + \left< \Nu; \vc{m} \right> \cdot \Grad \varphi \Big] } \dt
\end{equation}
for any $\varphi \in C^1([0,T] \times \Ov{\Omega})$ and a.a. $\tau \in (0,T)$;

\item
\begin{equation} \label{p2}
\begin{split}
&\left[ \intO{ \left< \Nu ; \vc{m} \right> \cdot \bfphi } \right]_{t = 0}^{t = \tau} \\&=
\int_0^\tau \intO{ \left[ \left< \Nu; \vc{m} \right> \cdot
\partial_t \bfphi + \left< \Nu; \frac{ \vc{m} \otimes \vc{m} }{\vr} \right> : \Grad \bfphi +
\frac{1}{\ep^2} \left< \Nu; p \right> \ \Div \bfphi \right] } \dt\\
&+ \frac{1}{\ep^2} \int_0^\tau \intO{ \left< \Nu ;\vc{m} \right> \cdot \Grad F } \dt  + \int_0^\tau \int_{\Ov{\Omega}} \Grad \bfphi : {\rm d} \mu_c
\end{split}
\end{equation}
for any $\bfphi \in C^1([0,T] \times \Ov{\Omega}, R^N)$, $\bfphi \cdot \vc{n}|_{\partial \Omega} = 0$, and a.a. $\tau \in (0,T)$, where $\mu_c \in \mathcal{M}([0,T] \times \Ov{\Omega}, R^N \times R^N)$ is the so--called \emph{momentum concentration measure};
\item
\begin{equation} \label{p3}
\left[ \intO{ \left< \Nu ; \frac{1}{2} \frac{ |\vc{m} |^2}{\vr} + \frac{1}{\ep^2} c_v p \right> } \right]_{t = 0}^{t = \tau} \leq
\frac{1}{\ep^2}
\int_0^\tau \intO{ \left< \Nu; \vc{m} \right>
\cdot \Grad F } \dt
\end{equation}
for a.a. $\tau \in [0,T]$;

\item
\begin{equation} \label{p4}
\begin{split}
&
\left[ \intO{ \left< \Nu ; \vr \chi (s(\vr, p)) \right> \varphi } \right]_{t = 0}^{t = \tau} \\ &\geq
\int_0^\tau \intO{ \Big[ \left< \Nu; \vr \chi (s(\vr, p)) \right> \partial_t \varphi +  \left< \Nu; \chi (s(\vr, p)) \vc{m} \right> \cdot \Grad \varphi \Big] } \dt
\end{split}
\end{equation}
for any $\varphi \in C^1([0,T] \times \Ov{\Omega})$, $\varphi \geq 0$, a.a. $\tau \in (0,T)$, and any
\[
\chi \in C(R), \ \chi \ \mbox{concave},\ \chi(S) \leq \Ov{\chi}, \ \mbox{for all} \ S \in R;
\]

\item
the concentration defect measure satisfies
\begin{equation} \label{p5}
\int_0^\tau \int_{\Ov{\Omega}} |{\rm d}\mu_c | \leq C \int_0^\tau \mathcal{D}(t) \ \dt \ \mbox{for a.a.}\ \tau \in (0,T),
\end{equation}
where
\[
\begin{split}
\mathcal{D}(\tau) &\equiv \intO{ \left[ \left< \mathcal{V}_{0,x} ; \frac{1}{2} \frac{ |\vc{m} |^2}{\vr} + \frac{1}{\ep^2} c_v p \right> -  \left< \mathcal{V}_{\tau,x} ; \frac{1}{2} \frac{ |\vc{m} |^2}{\vr} + \frac{1}{\ep^2} c_v p \right> \right]} \\&+ \frac{1}{\ep^2}
\int_0^\tau \intO{ \left< \Nu; \vc{m} \right>
\cdot \Grad F } \dt.
\end{split}
\]
In addition, we require
the constant $C$ in (\ref{p5}) to be independent of $\ep$.

\end{itemize}

The parameterized measure $\{ \mathcal{V}_{0,x} \}_{x \in \Omega}$ plays the role of the initial data and is given. The above definition is slightly more restrictive than its counterpart in \cite{BreFei17}. In particular, the renormalizing functions $\chi$ in the entropy inequality (\ref{p4}) need not be increasing. Such a concept fits better to problems
arising as asymptotic limits of isentropic viscous systems, cf. Kr{\"o}ner, and Zajaczkowski \cite{KrZa}. We refer to \cite{BreFei17B}, \cite{BreFei17} for the motivation and the basic properties of the (DMV) solutions to the complete Euler system.
We recall the definition of the non-linearities appearing in (\ref{p2}), (\ref{p3}), (\ref{p4}) on the singular set $\vr = 0$ and/or $p = 0$:
\begin{equation} \label{p4a}
[\vr, \vc{m}] \mapsto \frac{ |\vc{m}|^2 }{\vr} = \left\{ \begin{array}{lcl} \frac{ |\vc{m}|^2  }{\vr} \ &&\mbox{if} \ \vr > 0, \\
\infty \ &&\mbox{if}\ \vr = 0 \ \mbox{and}\ \vc{m} \ne 0,\\ 0 \ &&\mbox{otherwise},    \end{array} \right.
\end{equation}
\begin{equation} \label{p4b}
[\vr, p] \mapsto \vr \log \left( \frac{p^{c_v}}{\vr^{1 + c_v}} \right) =\vr c_v \log \left( \frac{p}{\vr^\gamma} \right) =
 \left\{ \begin{array}{lcl} \vr c_v \log \left( \frac{p}{\vr^\gamma } \right)
\ &&\mbox{if} \ \vr \geq 0,\ p > 0, \\
- \infty \ &&\mbox{if}\ \vr > 0,\ p = 0,\\
0 \ &&\mbox{otherwise.} \end{array} \right.
\end{equation}
In particular, we tacitly assume in (\ref{p3}), (\ref{p4}) that the Borel functions defined through (\ref{p4a}), (\ref{p4b}) are $\Nu$-integrable for a.a. $(t,x)$ which implies
\begin{equation} \label{p5a}
\Nu \left\{ [\vr, \vc{m}, p ] \ \Big|\ \vr = 0 \ \mbox{and}\ \vc{m} \ne 0 \right\} =
\Nu \left\{ \vr > 0
\ \mbox{and}\ p = 0 \right\} = 0\ \mbox{for a.a.}\ (t,x).
\end{equation}

\subsection{A comparison principle for the entropy}
\label{CP}

The fact that the entropy $s$ satisfies, formally, the transport equation
\[
\partial_t s + \vu \cdot \Grad s = 0
\]
remains encoded in the (DMV) formulation (\ref{p4}). Indeed suppose that
\[
s_* \leq
s(\vr(0, \cdot), p(0,\cdot)) \leq s^* \ \mbox{ a.e. in } \ \Omega.
\]
In terms of the probability measure $\mathcal{V}_{0,x}$ this means that
\begin{equation} \label{p6}
\mathcal{V}_{0,x} \left\{ [\vr , \vc{m}, p ] \ \Big| s_* \leq \log \left( \frac{p^{c_v}}{\vr^{1 + c_v}} \right) \leq s^* \right\} = 1
\ \mbox{ for a.a. }\ x \in \Omega.
\end{equation}

Now, consider a concave function $\chi$,
\[
\chi_\Lambda (s) = \left\{ \begin{array}{lcl} \Lambda (s - s_* ) \ &\mbox{if}&\ s \leq s_*, \\ 
0 \ &\mbox{if}&\ s_* \leq s \leq s^*, \\ 
- \Lambda (s - s^* ) \ &\mbox{if}&\ s \geq s^*, \end{array} \right.  \qquad \Lambda > 0.
\]
Plugging $\chi_\Lambda$ in (\ref{p4}) and using (\ref{p6}), we obtain
\[
\intO{ \left< \Nu ; \vr \chi_\Lambda \left( s(\vr, p ) \right) \right> } \geq 0 \ \mbox{for a.a.} \ \tau \in (0,T).
\]
Consequently, relation (\ref{p6}) gives rise to
\begin{equation} \label{p7}
\mathcal{V}_{t,x} \left\{ [\vr , \vc{m}, p ] \ \Big| s_* \leq \log \left( \frac{p^{c_v}}{\vr^{1 + c_v}} \right) \leq s^* \right\}
\Big| \Big\{ \vr > 0 \Big\} = 1 \quad \mbox{for a.a.}\ (t,x).
\end{equation}
This means that $s$ remains in the strip $[ s_*, s^* ]$ as long as $\vr > 0$.
In particular, if the initial data satisfy (\ref{p6}), we may use an approximation argument to relax the hypothesis $\chi \leq \Ov{\chi}$ in (\ref{p4}). Thus (\ref{p4}) holds for \emph{any} concave $\chi$; whence the choice $\chi(s) = \pm s$ gives rise
\begin{equation} \label{p8}
\begin{split}
&
\left[ \intO{ \left< \Nu ; \vr s(\vr, p) \right> \varphi} \right]_{t = 0}^{t = \tau} \\ &=
\int_0^\tau \intO{ \Big[ \left< \Nu; \vr s(\vr, p) \right> \partial_t \varphi +  \left< \Nu; s(\vr, p) \vc{m} \right> \cdot \Grad \varphi \Big] } \dt
\end{split}
\end{equation}
for any $\varphi \in C^1([0,T] \times \Ov{\Omega})$, and a.a. $\tau \in (0,T)$.

\begin{Remark} \label{RM1}

Note that (\ref{p7}) yields
\begin{equation} \label{p7a}
\Nu \left\{ [\vr, \vc{m}, p] \ \Big|\ \exp \left( \frac{s_*}{c_v} \right) \vr^\gamma \leq p \right\} = 1,
\end{equation}
in particular, in view of the energy bound (\ref{p3}),
the quantities $\left< \Nu ; \vr s(\vr, p) \right>$ and $\left< \Nu; s(\vr, p) \vc{m} \right>$ are integrable in $(0,T) \times \Omega$.

\end{Remark}

\section{Main result}
\label{M}

We formulate our main result concerning the singular limit of the system (\ref{i1}--\ref{i3}) for the initial density
and temperature distribution close to the isentropic static state. More specifically, we suppose that
\begin{equation} \label{ss1}
\Grad \tilde{p} = \tvr \Grad F \ \mbox{in}\ \Omega,\ \mbox{with}\ \tilde{p} = \exp \left( \frac{ \Ov{s} }{c_v} \right) \tvr^\gamma ,\
\tvr > 0 \ \mbox{in}\ \Ov{\Omega},
\end{equation}
where $\bar s$ is a constant, $\Omega \subset R^N$, $N = 2,3$ is a bounded regular domain and $F$ is a smooth potential.

\begin{Theorem}\label{T1}
Let $\Omega \subset R^N$, $N=2,3$ be a bounded domain with smooth boundary.
Let $\tvr$, $\tilde p$ be the static solution determined by \eqref{ss1} and let
\[
\vc{U}_0\in W^{k,2}(\Omega; R^N), k\geq N, \ \Div (\tvr \vc{U}_0) = 0, \ \vc{U}_0 \cdot \vc{n}|_{\partial \Omega} = 0.
\]
Suppose that
the anelastic Euler system (\ref{i13})--(\ref{i14}) with initial datum $\vc{U}_0$  admits a unique strong solution $\vU$ defined on a maximal time interval $[0, T_{\rm max})$.

Let $\{\Nu^\ep\}_{(t,x)\in (0,T)\times \Omega} $, $0 < T < T_{\rm max}$, be a family of dissipative measure-valued solutions to (\ref{i1}--\ref{i3}), satisfying \eqref{p5} with a constant $C$ independent of $\ep$, emanating from the initial data $\mathcal{V}_{0,x}^\ep$ so that
\begin{align*}
\mathcal{V}_{0,x}^\ep\left\{ [\vr , \vc{m}, p ] \ \Big| \left| \frac{\vr- \tvr}{\ep}\right| + \left|\frac{\vm}{\vr}- \vU_0\right|+ \left|\frac{p- \tilde{p} }{\ep}\right|\leq M_\ep(x) \right\} = 1
\ \mbox{ for a.a. }\ x \in \Omega,
\end{align*}
where
\[
 \|M_\ep\|_{L^\infty (\Omega)} \leq c \qquad \mbox{and}\qquad 	M_\ep \to 0 \mbox{ in } L^1(\Omega) \mbox{ as } \ep \to 0.
\]
Moreover, suppose that the initial entropy satisfies
\begin{equation}\label{t1}
\mathcal{V}_{0,x}^\ep\left\{ [\vr , \vc{m}, p ] \ \Big| \bar s- \ep^{2+\alpha}\leq s(\vr,p)\leq \bar s + \ep^{2+\alpha} \right\} = 1
\ \mbox{ for a.a. }\ x \in \Omega
\end{equation}
for some $\alpha>0$, where $\bar s = s(\tvr, \tilde{p})$ as in \eqref{ss1}.

Then
\[
\mathcal{D}^\ep \to 0 \quad \mbox{in } L^\infty(0,T) \ \mbox{as}\ \ep \to 0,
\]
and
\[
\mathcal{V}^\ep \to \delta_{\tilde{\vr},\tvr \vc{U}, \tilde{p} }\quad \mbox{in }L^\infty(0,T, L^q (\Omega; \mathcal{M}^+(\mathcal{F})_{\emph{weak}-^*}))
\ \mbox{as}\ \ep \to 0 \ \mbox{for any}\ 1 \leq q < \infty.
\]
\end{Theorem}

\begin{Remark} The regularity of the initial datum $\vU_0$ in Theorem \ref{T1} is motivated by the result by Oliver \cite{Oli}, where the existence of strong solutions for the anelastic Euler system is proved for sufficiently smooth domains $\Omega \subset R^2$ in the stated class of initial data.
\end{Remark}

\begin{Remark} \label{RM2}

The convergence
\[
\mathcal{V}^\ep \to \delta_{\tilde{\vr},\tvr \vc{U}, \tilde{p} }\quad \mbox{in }L^\infty(0,T, L^q(\Omega; \mathcal{M}^+(\mathcal{F})_{\emph{weak}-^*}))
\]
means that
\[
{\rm ess} \sup_{t \in (0,T)} \intO{ \left| \left< \Nu; G(\vr, \vc{m}, p) \right> - G(\tvr, \tvr \vc{U}, \tilde{p} )(t,x)
\right|^q } \to 0 \ \mbox{for any}\ G \in C_c (\mathcal{F}).
\]

\end{Remark}

\begin{Remark} \label{RM3}

In view of the comparison principle discussed in Section \ref{CP}, hypothesis (\ref{t1}) implies (\ref{p7}), (\ref{p7a}) with
\[
s_* = \Ov{s} - \ep^{2 + \alpha}, \ s^* = \Ov{s} + \ep^{2 + \alpha}.
\]
In particular, the entropy satisfies equation (\ref{p8}).

\end{Remark}

The proof of Theorem \ref{T1} is given in the next section.
Recall that if the space dimension is $N=2$, the smooth solution $\vU$ of the target system \eqref{i13}-\eqref{i14} exists globally in time and $T_{\rm max}=\infty$ in the above theorem.

\section{Proof of the main result}
\label{PM}

Similarly to \cite[Section 2.4.1]{BreFei17B}, we introduce the total entropy
\[
\mathcal{S}(\vr , p) = \vr \log \left( \frac{p^{c_v}}{ \vr^{1 + c_v} } \right) = c_v \vr \log \left( \frac{p}{ \vr^{\gamma} } \right), \ \gamma = 1 + \frac{1}{c_v},
\]
together with the \emph{relative energy}
\[
\begin{split}
\mathcal{E}_\ep &\left( \vr, \vc{m}, p \Big| r, \tvu, \Theta \right) \equiv \\
&\frac{1}{2} \frac{ |\vc{m} |^2 }{\vr} + \frac{1}{\ep^2} c_v p - \Theta
\frac{1}{\ep^2} \mathcal{S}(\vr, p ) - \vc{m} \cdot \tvu + \frac{1}{2} \vr |\tvu|^2 + \frac{1}{\ep^2} P(r, \Theta) - \frac{\vr}{\ep^2} \left( E(r, \Theta) - \Theta S(r, \Theta) + \frac{P(r, \Theta)}{r} \right).
\end{split}
\]
 The
thermodynamic functions $E$, $P$, and $S$ - the internal energy, the pressure, and the entropy - are expressed in terms of the standard
state variables,
\[
E (r, \Theta) = c_v \Theta, \ P(r, \Theta) = r \Theta, \ S(r, \Theta) = \log \left( \frac{\Theta^{c_v}}{r} \right).
\]

As shown in \cite[Section 2.4.1]{BreFei17B}, \cite{FeKlKrMa}, any DMV solution to the Euler system satisfies the relative energy inequality
\begin{equation} \label{p9}
\begin{split}
\Big[ \int_\Omega &\left< \Nu^\ep ; \mathcal{E}_\ep \left( \vr, \vc{m}, p \Big| r, \tilde \vu, \Theta \right) \right> \ \dx \Big]_{t = 0}^{t = \tau}
+ \mathcal{D}^\ep (\tau) \\
\leq &- \frac{1}{\ep^2} \int_0^\tau \intO{ \left[ \left< \Nu^\ep ; \vr s(\vr, p) \right> \partial_t \Theta +
\left< \Nu^\ep ;  s(\vr, p) \vc{m} \right> \cdot \Grad \Theta \right] } \dt\\
&+ \int_0^\tau \intO{ \left[ \left< \Nu^\ep ; \vr \tvu - \vc{m} \right> \cdot \partial_t \tilde \vu  + \left< \Nu^\ep ; \frac{ (\vr \tvu - \vc{m}) \otimes \vc{m} }{\vr}
\right> : \Grad \tilde \vu  \right] } \dt \\ &- \frac{1}{\ep^2} \int_0^\tau \intO{ \left[ \left< \Nu^\ep ; p \right> \Div \tilde \vu \right] }\dt \\
&+ \frac{1}{\ep^2} \int_0^\tau \intO{ \left[ \left< \Nu^\ep ; \vr \right> \partial_t \Theta S(r, \Theta) + \left< \Nu^\ep ; \vc{m} \right> \cdot \Grad \Theta S(r, \Theta) \right] } \dt\\
&+ \frac{1}{\ep^2} \int_0^\tau \intO{ \left[ \left< \Nu^\ep ; r - \vr \right> \frac{1}{r} \partial_t P (r, \Theta) - \left< \Nu^\ep ; \vc{m} \right> \cdot \frac{1}{r} \Grad P(r, \Theta) \right] } \dt\\
&+ \frac{1}{\ep^2} \int_0^\tau \intO{ \Grad F \cdot \left< \Nu^\ep; \vc{m} - \vr \tvu \right> } \dt +  \int_0^\tau \int_{\bar \Omega} \Grad \tvu : {\rm d} \mu_{c,\ep}
\end{split}
\end{equation}
as soon as the initial data satisfy (\ref{p6}).
Relation (\ref{p9}) holds for \emph{any} trio of differentiable test functions $[r, \tvu, \Theta]$ satisfying
\[
r > 0, \ \Theta > 0, \ \tvu \cdot \vc{n}|_{\partial \Omega} = 0
\]
that play the role of standard state variables - the density, the velocity, and the absolute temperature.

The proof of Theorem \ref{T1} will be done with the help of (\ref{p9}), where we take
\[
r = \tvr , \ \Theta = \tvt \equiv \frac{\tilde{p}}{\tvr},
\]
$\tvr$, $\tilde{p}$ being the static solution (\ref{ss1}). Accordingly, relation (\ref{p9}) simplifies to
\begin{equation} \label{p9a}
\begin{split}
\Big[ \int_\Omega &\left< \Nu^\ep ; \mathcal{E}_\ep \left( \vr, \vc{m}, p \Big| \tvr, \tilde \vu, \tvt \right) \right> \ \dx \Big]_{t = 0}^{t = \tau}
+ \mathcal{D}^\ep (\tau) \\ \leq &
\int_0^\tau \intO{ \left< \Nu^\ep ; \frac{ (\vr \tvu - \vc{m}) \otimes (\vc{m} - \vr \tvu)  }{\vr}
\right> : \Grad \tilde \vu } \dt
\\
&+ \frac{1}{\ep^2} \int_0^\tau \intO{
\left< \Nu^\ep ; \left(\Ov{s} -   s(\vr, p) \right) \vc{m} \right> \cdot \Grad \tvt } \dt\\
&+ \int_0^\tau \intO{ \left< \Nu^\ep ; \vr \tvu - \vc{m} \right> \cdot \left( \partial_t \tilde \vu +
\tvu \cdot \Grad \tvu \right) } \dt \\ &- \frac{1}{\ep^2} \int_0^\tau \intO{ \left< \Nu^\ep ; p \right> \Div \tilde \vu  }\dt \\
&- \frac{1}{\ep^2} \int_0^\tau \intO{ \Grad F \cdot \left< \Nu^\ep; \vr \tvu \right> } \dt +  \int_0^\tau \int_{\bar \Omega} \Grad \tvu : {\rm d} \mu_{c,\ep}.
\end{split}
\end{equation}
The proof of Theorem \ref{T1} will be carried over in the next two sections by taking
$\tvu = \vc{U}$ in (\ref{p9a}) and performing the limit $\ep \to 0$.
Notice that the initial datum $\mathcal{V}^\ep_0$ in Theorem \ref{T1} is chosen in a way such that
\begin{align}\label{p10}
	\int_\Omega &\left< \mathcal{V}_{0,x}^\ep ; \mathcal{E}_\ep \left( \vr, \vc{m}, p \Big| \tvr, \tvt, \vU_0 \right) \right> \ \dx \to 0 \quad \mbox{ as } \ep \to 0.
\end{align}

\subsection{Coerciveness and uniform estimates}

We start with deriving uniform bounds independent of $\ep \to 0$.

\subsubsection{Coerciveness of the relative energy}

As in \cite[Section 3.2]{BreFei17} we introduce the \emph{essential} and \emph{residual} part of a function $G(\vr, \vm,p)$, which describe the behavior of $G$ in the area where $\vr$ and $p$ are bounded from below and above, and outside of this area, respectively.
Let $[\tvr, \tilde{p}]$ be the solution of (\ref{ss1}), and $\tilde K$ be its image in $(0,\infty)^2$,
\[
\tilde{K} = \left\{ [\tvr, \tilde{p}] (x) \ \Big| \ x \in \Ov{\Omega} \right\}.
\]
Next, we introduce a cutoff function
\[
\Psi \in C_c^\infty ((0,\infty)^2), \quad 0\leq \Psi \leq 1, \quad \Psi |_\mathcal{U}=1, \mbox{ where } \mathcal{U}\subset (0,\infty)^2 \mbox{ is an open set containing } \tilde K.
\]
Now, any (Borel) function $G(\vr, \vm, p)$ can be decomposed as a sum of its \emph{essential} and \emph{residual} part as
\[
G(\vr, \vm, p)= [G(\vr, \vm, p)]_{ess} + [G(\vr, \vm, p)]_{res},
\]
where
\[
[G(\vr, \vm, p)]_{ess}:= \Psi(\vr, p)G(\vr, \vm, p) \qquad \mbox{and}\qquad [G(\vr, \vm, p)]_{res}:=(1-\Psi(\vr, p))G(\vr, \vm, p).	
\]
As in \cite[Section 3.2.2]{BreFei17}, the relative energy is a coercive functional satisfying the estimate
\begin{align}\label{p11}
\begin{split}
&\mathcal{E}_\ep \left( \vr, \vc{m}, p \Big| \tilde \vr, \tvu , \tvt  \right) \gtrsim \left[\left|\frac{\vm}{\vr}-\tvu \right|^2\right]_{ess}+ \left[\frac{|\vm|^2}{\vr} \right]_{res}\\
&\qquad+ \frac{1}{\ep^2}\left[ |\vr-\tvr|^2 + | p - \tilde{p} |^2 \right]_{ess}+ \frac{1}{\ep^2}\left[ 1+ \vr + \vr|s(\vr,p)|+ p \right]_{res}, \
\tilde{p} = P(\tvr, \tvt) = \tvr \tvt.
\end{split}
\end{align}

\subsubsection{Uniform bounds}

Taking $\tvu=0$ in \eqref{p9a} we obtain
\begin{equation*}
\begin{split}
\Big[ \int_\Omega &\left< \Nu^\ep ; \mathcal{E}_\ep \left( \vr, \vc{m}, p \Big| \tilde \vr, 0, \tvt \right) \right> \ \dx \Big]_{t = 0}^{t = \tau}
+ \mathcal{D}^\ep (\tau) \\
\leq & \frac{1}{\ep^2} \int_0^\tau \intO{
\left< \Nu^\ep ;  \vc{m}\left( \Ov{s}-s(\vr, p) \right) \right> \cdot \Grad \tilde \vt  } \dt.
\end{split}
\end{equation*}
In view of \eqref{p7}, the entropy inequality \eqref{p4} together with \eqref{t1}  guarantees that
\begin{align*}
  \intO{
\left< \Nu^\ep ; \vr \left(\frac{\bar s -s(\vr, p)}{\ep^2}\right)^2 \right>   }  \leq \intO{
\left< \mathcal{V}_{0,x}^\ep ; \vr \left(\frac{\bar s -s(\vr, p)}{\ep^2}\right)^2 \right>   } \lesssim 1.
\end{align*}
Therefore,
\begin{align*}
\Big[ \int_\Omega \left< \Nu^\ep ; \mathcal{E}_\ep \left( \vr, \vc{m}, p \Big| \tilde \vr, \tilde \vt, 0 \right) \right> \ \dx \Big]_{t = 0}^{t = \tau}
+ \mathcal{D}^\ep (\tau)
\lesssim   \int_0^\tau \intO{
\left< \Nu^\ep ;  \frac{1}{2}\frac{|\vc{m}|^2}{\vr}\right>  } \dt + 1.
\end{align*}
Taking \eqref{p11} into account, an application of Gronwall's lemma provides that
\begin{equation*}
 \int_\Omega\left< \mathcal{V}_{\tau,x}^\ep ; \mathcal{E}_\ep \left( \vr, \vc{m}, p \Big| \tilde \vr, 0, \tvt \right) \right> \ \dx
+ \mathcal{D}^\ep (\tau)
\leq c.
\end{equation*}
where $c>0$ is a constant independent of $\ep$.
The above estimate gives rise to the following bounds:
\begin{align}\label{p17}
	\int_\Omega\left< \mathcal{V}_{\tau,x}^\ep ; \frac{|\vm|^2}{\vr} \right> \ \dx \leq c
\end{align}
and
\begin{align}\label{p12}
	\int_\Omega\left< \mathcal{V}_{\tau,x}^\ep ; [\vr-\tvr]^2_{ess} + [p- \tilde{p} ]^2_{ess} + 1 + \vr_{res}
+ p_{res} \right> \ \dx \lesssim \ep^2c \ \mbox{uniformly for}\ \tau \in [0,T].
\end{align}
Eventually, due to the same argument as in \cite[Section 4.4.3]{FeKlKrMa}, there exists a function $\vM \in L^\infty(0,T;L^q(\Omega))$ for some $q>1$ and a subsequence (not relabeled) such that
\begin{equation}\label{p13}
	\left< \Nu^\ep ; \vm \right> \to \vM\quad \mbox{weakly-* in } L^\infty(0,T;L^q(\Omega)).
\end{equation}
In view of the convergence \eqref{p13} and the uniform estimate \eqref{p12}, the conservation of mass equality \eqref{p1} ensures that
\[
	\int_0^\tau \int_{\Omega} \vM \cdot \nabla_x \varphi \dx \dt =0
\]
for any $\varphi \in C^1([0,T]\times \overline{\Omega})$ and a.a. $\tau \in (0,T)$. In other words
\[
	\Div \vM=0, \qquad \vM \cdot \vn \big|_{\partial \Omega}=0
\]
in the sense of distributions.

\subsection{The limit for $\ep\to 0$}
Our aim is to show that
\begin{equation*}
 \lim_{\ep \to 0} \left(\int_\Omega \left< \mathcal{V}_{\tau,x} ; \mathcal{E}_\ep \left( \vr, \vc{m}, p \Big| \tilde \vr, \vU, \tilde \vt \right) \right> \ \dx
+ \mathcal{D}^\varepsilon (\tau) \right) =0
\end{equation*}
uniformly for a.a. $\tau \in (0,T)$, where $\tvr, \tvt= \tilde{p} / \tvr$ are the isentropic static state solutions and $\vU$ the corresponding unique solution to the anelastic Euler system \eqref{i13}-\eqref{i14}. Then the coervivity property \eqref{p11} of $\mathcal{E}$ allows us to conclude the proof of Theorem \ref{T1}.

Let us denote by $\omega(\ep)$ a generic function, which may differ from occurrence to occurrence enclosing terms that vanish in the limit $\ep \to 0$.
For $\tvu = \vU$,
the relative energy inequality \eqref{p9a} reads
\[
\begin{split}
\int_\Omega &\left< \mathcal{V}_{\tau,x} ^\ep ; \mathcal{E}_\ep \left( \vr, \vc{m}, p \Big| \tvr, \vU , \tvt \right) \right>  \dx
+ \mathcal{D}^\ep (\tau) \\ &\leq
\int_0^\tau \intO{ \left< \Nu^\ep ; \frac{ (\vr \vU - \vc{m}) \otimes (\vc{m} - \vr \vU)  }{\vr}
\right> : \Grad \vU } \dt
\\
&\quad+ \frac{1}{\ep^2} \int_0^\tau \intO{
\left< \Nu^\ep ; \left(\bar{s} -   s(\vr, p) \right) \vc{m} \right> \cdot \Grad \tvt } \dt\\
&\quad- \int_0^\tau \intO{ \left< \Nu^\ep ; \vr \vU - \vc{m} \right> \cdot \Grad \Pi } \dt \\ &\quad - \frac{1}{\ep^2} \int_0^\tau \intO{ \left< \Nu^\ep ; p \right> \Div \vU  }\dt \\
&\quad- \frac{1}{\ep^2} \int_0^\tau \intO{ \Grad F \cdot \left< \Nu^\ep; \vr \vU \right> } \dt +  \int_0^\tau \int_{\bar \Omega} \Grad \vU : {\rm d} \mu_{c,\ep}
+ \omega(\ep),
\end{split}
\]
where we have used \eqref{p10}.

Furthermore, in view of \eqref{p5}, (\ref{p11}),
\[
\begin{split}
\int_0^\tau &\intO{ \left< \Nu^\ep ; \frac{ (\vr \vU - \vc{m}) \otimes (\vc{m} - \vr \vU)  }{\vr}
\right> : \Grad \vU } \dt + \int_0^\tau \int_{\bar \Omega} \Grad \vU : {\rm d} \mu_{c,\ep}\\
&\lesssim \int_0^\tau \left[ \int_\Omega \left< \Nu^\ep ; \mathcal{E}_\ep \left( \vr, \vc{m}, p \Big| \tvr, \vU, \tvt \right) \right>  \dx
+ \mathcal{D}^\ep(t) \right] \dt;
\end{split}
\]
whence
\begin{equation} \label{p9b}
\begin{split}
\int_\Omega &\left< \mathcal{V}_{\tau,x} ^\ep ; \mathcal{E}_\ep \left( \vr, \vc{m}, p \Big| \tvr, \vU, \tvt \right) \right>  \dx
+ \mathcal{D}^\ep (\tau) \\ &\lesssim
\frac{1}{\ep^2} \int_0^\tau \intO{
\left< \Nu^\ep ; \left(\Ov{s} -   s(\vr, p) \right) \vc{m} \right> \cdot \Grad \tvt } \dt\\
&\quad - \int_0^\tau \intO{ \left< \Nu^\ep ; \vr \vU - \vc{m} \right> \cdot \Grad \Pi } \dt \\ &\quad- \frac{1}{\ep^2} \int_0^\tau \intO{ \left< \Nu^\ep ; p \right> \Div \vU  }\dt - \frac{1}{\ep^2} \int_0^\tau \intO{ \Grad F \cdot \left< \Nu^\ep; \vr \vU \right> } \dt\\
&\quad+ \int_0^\tau \left[ \int_\Omega \left< \Nu^\ep ; \mathcal{E}_\ep \left( \vr, \vc{m}, p \Big| \tvr, \vU, \tvt \right) \right>  \dx
+ \mathcal{D}^\ep(t) \right] \dt
+ \omega(\ep).
\end{split}
\end{equation}

We are going to consider the remaining integrals on the right hand side of \eqref{p9b} step by step.

\subsubsection*{Step 1}
Let $\bar \alpha \in (0,\alpha)$, where $\alpha>0$ is the constant in hypothesis \eqref{t1}. Then
\begin{align}\label{p14}
\begin{split}
 \frac{1}{\ep^2} \int_0^\tau& \intO{
\left< \Nu^\ep ;  \left(\bar s-s(\vr, p) \right) \vc{m} \right> \cdot \Grad \tilde \vt  } \dt\\
&\lesssim   \int_0^\tau \intO{
\left< \Nu^\ep ;  \left( \frac{\bar s-s(\vr, p)}{\ep^{2+\bar \alpha}} \right)^2 \vr \right>   } \dt +  \ep^{2 \bar \alpha} \int_0^\tau \intO{
\left< \Nu^\ep ;  \frac{|\vc{m}|^2}{\vr} \right>   } \dt\\
&=\omega(\ep),
\end{split}
\end{align}
where we used \eqref{p4}, \eqref{p7}, \eqref{t1} for the entropy integral, and the uniform bound \eqref{p17}.

\subsubsection*{Step 2}

Next, in accordance with (\ref{i13}),
\[
\int_0^\tau \intO{ \left< \Nu; \vr \vU - \vc{m} \right> \Grad \Pi } \dt = \int_0^\tau \intO{ \left< \Nu; (\vr - \tvr) \vU - \vc{m} \right> \Grad \Pi } \dt.
\]
Consequently, in view of (\ref{p12}), (\ref{p13}),
\begin{equation} \label{p14a}
\int_0^\tau \intO{ \left< \Nu; \vr \vU - \vc{m} \right> \Grad \Pi } \dt \lesssim
\int_0^\tau \int_\Omega \left< \Nu^\ep ; \mathcal{E}_\ep \left( \vr, \vc{m}, p \Big| \tvr, \vU, \tvt \right) \right>  \dx \ \dt
+ \omega(\ep).
\end{equation}

\subsubsection*{Step 3}
The relations $\nabla_x F = \frac{1}{\tvr}\nabla_x \tilde p$ and $\vU \cdot \vc{n} \big|_{\partial \Omega}=0$ give rise to
\begin{align}\label{eq:add}
\begin{split}
&-\frac{1}{\ep^2}\int_0^\tau \intO{  \left< \Nu^\ep ; p \right> \Div  \vU  }\dt-\frac{1}{\ep^2}\int_0^\tau \intO{ \Grad F \cdot \left< \Nu^\ep;  \vr \vU \right> } \dt\\
&=-\frac{1}{\ep^2}\int_0^\tau \intO{  \left< \Nu^\ep ; p - \tilde{p}  \right> \Div  \vU  }\dt \\
&\quad + \frac{1}{\ep^2}\int_0^\tau \intO{ \left[\nabla_x \tilde{p} \cdot \vU - \frac{1}{\tvr}\nabla_x \tilde{p} \cdot \left< \Nu^\ep;  \vr \vU \right> \right] } \dt\\
&= -\frac{1}{\ep^2}\int_0^\tau \intO{  \left< \Nu^\ep ; p - \tilde p  \right> \Div  \vU  }\dt -\frac{1}{\ep^2}\int_0^\tau \intO{ \left< \Nu^\ep;  \vr-\tvr  \right> \frac{1}{\tvr}\nabla_x \tilde p \cdot \vU  } \dt.
\end{split}
\end{align}
Notice that the residual part of the above terms can be controlled by $ \mathcal{E}_\ep \left( \vr, \vc{m}, p \Big| \tilde \vr, \vU, \tilde \vt \right)$.

Concerning the essential part, we use hypothesis \eqref{t1}, together with the corresponding relations (\ref{p7}), (\ref{p7a}) to obtain
\[
\begin{split}
-\frac{1}{\ep^2}&\int_0^\tau \intO{  \left< \Nu^\ep ; [ p - \tilde p ]_{ess}  \right> \Div  \vU  }\dt -\frac{1}{\ep^2}\int_0^\tau \intO{ \left< \Nu^\ep;
[\vr-\tvr]_{ess}  \right> \frac{1}{\tvr}\nabla_x \tilde p \cdot \vU  } \dt \\&= 
-\frac{1}{\ep^2}e^{\frac{\bar s}{c_v}}\int_0^\tau \intO{  \left< \Nu^\ep ; [ \vr^\gamma - \tvr^\gamma ]_{ess}  \right> \Div  \vU  }\dt \\ &\quad -\frac{1}{\ep^2}e^{\frac{\bar s}{c_v}}\int_0^\tau \intO{ \left< \Nu^\ep;  [ \vr-\tvr ]_{\rm ess}  \right> \gamma \tvr^{\gamma - 2} \Grad \tvr \cdot \vU  } \dt + \omega(\ep)\\ &=
-\frac{1}{\ep^2}e^{\frac{\bar s}{c_v}}\int_0^\tau \intO{  \left< \Nu^\ep ; [ \vr^\gamma - \gamma \tvr^{\gamma - 1} (\vr - \tvr) - \tvr^\gamma ]_{ess}  \right> \Div  \vU  }\dt
+ \omega(\ep),
\end{split}
\]
where we have used the anelastic constraint (\ref{i13}) and
\[
	[p-\tilde p]_{ess} = e^{\frac{\bar s}{c_v}}[\vr^\gamma - \tvr^\gamma]_{ess}+ \frac{\tilde p}{c_v}[s(\vr,p)- \bar s]_{ess}+c\left( [\vr^\gamma - \tvr^\gamma]^2_{ess} + [p-\tilde p]^2_{ess} \right).
\]
Thus we may infer
\begin{equation} \label{p15}
\begin{split}
-\frac{1}{\ep^2}&\int_0^\tau \intO{  \left< \Nu^\ep ; [ p - \tilde p ]_{ess}  \right> \Div  \vU  }\dt -\frac{1}{\ep^2}\int_0^\tau \intO{ \left< \Nu^\ep;
[\vr-\tvr]_{ess}  \right> \frac{1}{\tvr}\nabla_x \tilde p \cdot \vU  } \dt \\&\lesssim
\int_0^\tau \int_\Omega \left< \Nu^\ep ; \mathcal{E}_\ep \left( \vr, \vc{m}, p \Big| \tvr, \vU, \tvt \right) \right>  \dx \ \dt
+ \omega(\ep).
\end{split}
\end{equation}

Collecting \eqref{p14}-\eqref{p15}, we may conclude that
\begin{equation*}
\begin{split}
 \int_\Omega &\left< \mathcal{V}_{\tau,x}^\ep ; \mathcal{E}_\ep \left( \vr, \vc{m}, p \Big| \tilde \vr, \vU, \tvt \right) \right> \ \dx
+ \mathcal{D}^\ep (\tau) \\
&\qquad \lesssim \int_0^\tau \intO{\left< \Nu^\ep ;\mathcal{E}_\ep \left( \vr, \vc{m}, p \Big| \tilde \vr, \vU, \tvt \right) \right> } + \mathcal{D}^\ep(t) \dt + \omega(\ep).
\end{split}
\end{equation*}

Applying Gronwall's lemma yields the desired estimate
\begin{align*}
\intO{\left< \Nu^\ep ;\mathcal{E}_\ep \left( \vr, \vc{m}, p \Big| \tilde \vr, \vU, \tvt \right) \right> } + \mathcal{D}^\ep(t)  \leq \omega(\ep)\quad \to 0 \qquad \mbox{as } \ep \to 0.
\end{align*}

We have proved Theorem \ref{T1}.

\section{Applications to the shallow water equation and related problems}
\label{A}

The \emph{shallow water equations} with bottom topography read:
\begin{equation}
\label{A1}
\begin{split}
\partial_t h + \Div (h \vu) &= 0,\\
\partial_t (h \vu) + \Div (h \vu \otimes \vu) + \frac{1}{\ep^2} h \Grad h &= \frac{1}{\ep^2} h \Grad b,
\end{split}
\end{equation}
where $\vu=(u_1,u_2)$ denotes the velocity field, $h$ is the fluid height, and $b$ represents the bottom topography, which is a given function depending on the space variables, see Lannes \cite{Lann} or Pedlosky \cite{PEDL}. System
(\ref{A1}) may be supplemented by the impermeability boundary condition
\begin{equation} \label{A2}
\vu \cdot \vc{n}|_{\partial \Omega} = 0,
\end{equation}
where $\Omega \subset R^2$ is a bounded domain. The singular limit $\ep \to 0$ represents the regime in which the gravitational force is dominant.

The system of equations \eqref{A1}, \eqref{A2} may be seen as a special case of (\ref{i1}--\ref{i3}) with constant entropy $s(\vr, \vt) = \Ov{s}$, $\gamma = 2$, and the obvious identification
$h \approx \vr$, $b \approx F$. Note that, in view of (\ref{p7}), the condition $s(\vr, \vt) = \Ov{s}$ may be enforced through a suitable choice of the initial data. Accordingly, the conclusion of
Theorem \ref{T1} applies to (\ref{A1}) as well.

Similar results may be obtained for certain modifications of (\ref{A1}) including the effect of blowing wind, see e.g. Csanady \cite{Csan}.

\def\cprime{$'$} \def\ocirc#1{\ifmmode\setbox0=\hbox{$#1$}\dimen0=\ht0
  \advance\dimen0 by1pt\rlap{\hbox to\wd0{\hss\raise\dimen0
  \hbox{\hskip.2em$\scriptscriptstyle\circ$}\hss}}#1\else {\accent"17 #1}\fi}

%\bibliographystyle{plain}
%\bibliography{citace}

\end{document}